  \documentclass{gretsi}
 \usepackage[latin1]{inputenc}   %
 \usepackage[english,french]{babel}   %
  \usepackage{times}			%

\usepackage{graphicx}
\usepackage{color}       %

\usepackage{verbatim}
\usepackage{amsmath}
\usepackage{amsfonts}
\usepackage{url}
\def\NLM{\mbox{NL-means}}
\def\s2{\sigma^2}
\def\Oxf{\Omega_{x}^{f}}
\def\Oxfronde{\Omega_{x}^{f^\circ}}
\def\Oxup{\Omega_{x}^{u^\circ}}
\def\EQM{\text{EQM}}
\def\Biais{\text{Biais}}
\def\Var{\text{Variance}}
\def\Covar{\text{Covariance}}
\newcommand\R{\mathbb R}

\titre{{\em Non-Local means} est un algorithme de débruitage local}

\auteur{\coord{Simon}{Postec}{},
        \coord{Jacques}{Froment}{},
    \coord{Béatrice}{Vedel}{}}

\adresse{
\affil{}{Université de Bretagne-Sud, UMR CNRS 6205,
Laboratoire de Mathématiques de Bretagne Atlantique, \\ Campus de Tohannic BP 573, F-56017 Vannes}
}

\email{Prenom.Nom@univ-ubs.fr}

\resumefrancais{L'algorithme de débruitage par moyennes non locales {\em Non-Local Means} ($\NLM$) introduit en 2005 \cite{bcm2005}
a repoussé les limites du débruitage mais il a aussi introduit un nouveau paradigme, selon lequel la similarité présente 
dans toute image naturelle peut être captée à travers les poids $\NLM$. Nous montrons que, contrairement à 
l'opinion dominante \cite{bcm2010}, les poids $\NLM$ ne permettent pas d'obtenir dans du bruit une mesure fiable de la similarité sans 
ajout d'une contrainte de localité : $\NLM$ considéré comme une méthode de débruitage s'avère être un algorithme local. Si quelques travaux \cite{go2007, gzw2011, s2010, t2009} avaient déjà remarqué que les meilleures performances de $\NLM$ s'obtenaient localement, aucune conclusion générale n'avait encore été proposée faute d'étude spécifique et la seule explication évoquée pour justifier les résultats expérimentaux s'avère insuffisante. Notre étude basée sur des expérimentations établit qu'en moyenne sur des images naturelles, le biais de l'estimateur $\NLM$ croît avec le rayon de la zone de recherche, ceci parce que le bruit perturbe l'ordre de similarité entre les patchs. 
Ainsi, l'erreur quadratique entre l'image originale et l'estimation $\NLM$ qui s'écrit comme la somme du biais, de la variance et de la covariance de l'estimateur possède un minimum absolu pour un disque de rayon 3 à 4 pixels.}

\resumeanglais{{\em Non-Local means} is a local image denoising algorithm. The {\em Non-Local Means} ($\NLM$) image denoising algorithm (\cite{bcm2005}, 2005) pushed the limits of denoising. But it introduced a new paradigm, according to which one could capture the similarity of images with the $\NLM$ weights. We show that, contrary to the prevailing opinion, the $\NLM$ weights do not allow to get a reliable measure of the similarity in a noisy image, unless one add a locality constraint. As an image denoising method, the {\em Non-Local Means} prove to be local. Some works \cite{go2007, gzw2011, s2010, t2009} had already pointed out that to get the best denoising performances with the $\NLM$ algorithm, one should run it locally. But no general conclusion has been yet proposed and the only explanation that was proposed to justify the experimental results is not sufficient. Our study based on experimental evidence proves that, on average on natural images, the bias of the $\NLM$ estimator is an increasing function of the radius of the similarity searching zone. The reason for this phenomenon is that noise disrupts the order of similarity between patches. Hence the mean squared error between the original image and the $\NLM$ estimation, which is the sum of the bias, the variance and the covariance of the estimator, has an absolute minimum for a disk of radius 3 to 4 pixels.}

\begin{document}
\maketitle

\section{L'algorithme $\NLM$}

Soient $v : \Omega \longrightarrow \R$ l'image bruitée et $u$ l'image originale. Nous supposons que $v = u + b$, où $\{b(x)\}_{x \ \in \ \Omega}$ 
est un bruit blanc gaussien de variance $\s2$. Dans les expériences qui suivent $\sigma=20$, mais d'autres valeurs usuelles de $\sigma$ 
conduisent aux mêmes conclusions qualitatives. Le débruitage par l'algorithme $\NLM$ \cite{bcm2005,bcm2010} est l'estimateur
\begin{equation}
\label{eq:nlm_filter}
\b{u}(x) = \displaystyle{\sum_{y \ \in \ \Omega_{x} } w(x,y) v(y)}
\end{equation}
où les poids $\NLM$ s'écrivent 
\begin{equation}
\label{eq:nlm_weights}
w(x,y) = \frac{1}{Z(x)} e^{-\frac{\|V(x) - V(y)\|_{2,a}^{2}}{2{h}^2}}.
\end{equation}
Nous notons $F(x)$, où $F$ peut être $U$, $V$ ou $B$, le patch de taille $7 \times 7$ centré en $x$ dans l'image $f=u,v$ ou $b$, c'est-à-dire les niveaux de gris $\{f(y) \ | \ y\in\Omega,{\| x - y \|}_\infty \leq 3\}$; $h$ est un paramètre de filtrage ici égal à $\sigma$, $Z(x)$ est un coefficient de normalisation, ${\| \ . \ \|}_{2,a}$ est la norme euclidienne pondérée par une gaussienne d'écart-type $a$ 
définissant la similarité entre les patchs et $\Omega_x = \{ \ y  \in  \Omega \ | \ {\| x - y \|}_2 \leq d \}$ est le voisinage du point 
$x$ de rayon $d$, correspondant à la zone de recherche de similarités.
Précisons que les articles originaux \cite{bcm2005,bcm2010} considèrent $\Omega_x=\Omega$ pour tout $x \ \in \ \Omega$; 
le choix de se restreindre à un voisinage $\Omega_x$ est conseillé au seul motif de limiter le temps de calcul.

\section{Variantes {\em a priori } mieux adaptées à la non-localité}

Des travaux \cite{go2007,gzw2011,s2010,t2009} ont mentionné, sur quelques images, la décroissance du rapport signal à bruit de crête $\mbox{PSNR}(u, \b{u})$ 
entre l'image bruitée $u$ et son estimateur $\b{u}$ pour des zones de recherche étendues et des versions dites semi-locales de l'algorithme ont été préférées.
Seuls \cite{s2010} et \cite{dag2011} proposent une explication: l'accumulation, dans la somme (\ref{eq:nlm_filter}), de contributions $v(y)$ affectées 
d'un petit poids. Un moyen simple d'éliminer ce phénomène est de mettre à zéro les poids inférieurs à un certain seuil ou de ne retenir que les premiers
poids, classés par ordre décroissant. Les patchs $V(y)$ pour lesquels $w(x,y)$ est retenu (ou est non nul) sont alors dits similaires à $V(x)$.
Il est intéressant de comparer ces patchs similaires avec ceux obtenus à partir de l'image originale $u$ (oracle).
Pour représenter ces différents poids nous notons $w_{v}$ (respectivement $w_{u}$) la fonction poids restreinte aux patchs similaires $V$
(respectivement aux patchs similaires oracle $U$). Ainsi, pour $f=u$ ou pour $f=v$ et pour $\Oxf$ l'ensemble $\Omega_x$ restreint
aux pixels $y$ tels que les patchs $F(y)$ et $F(x)$ soient similaires,
\begin{equation}
\label{eq:nlm_bestweights}
w_{f}(x,y) = \left\{ \begin{array}{ll}
w(x,y) & \mbox{ si } y \in \Oxf, \\
0 & \mbox{ sinon.}
\end{array}
\right.
\end{equation}
La notation $f^\circ$ à la place de $f$ signifie que le calcul de la distance $\|F(x) - F(y)\|_{2,a}$ ne tient pas compte de la valeur centrale des patchs, 
$f(x)$ et $f(y)$.
La suppression de la valeur centrale se justifie, dans le cas $f=v$, par le fait que ne conserver que les meilleurs patchs augmente le biais porté 
par le pixel central \cite{xxw2008}. 
Précisons que la fonction poids $w$ est toujours calculée selon l'équation (\ref{eq:nlm_weights}), donc à partir des patchs bruités $V$, le
cas échéant sans les valeurs centrales (notation $f^\circ$).
Dans ce qui suit les patchs similaires correspondent aux $80$ meilleurs patchs, mais les mêmes résultats qualitatifs sont obtenus avec un
seuillage. Le nombre $80$ est relativement arbitraire, d'autres valeurs donnent des résultats semblables.
Notons que quand $d \leq 5$, le voisinage  $\Omega_x$ contient moins de $80$ pixels et tous les patchs sont alors considérés comme similaires.

\section{Simulation avec les différents poids}

\begin{figure}[h]
\begin{center}
\includegraphics[width=80mm,bb=4cm 1.5cm 26cm 19cm,clip=true]{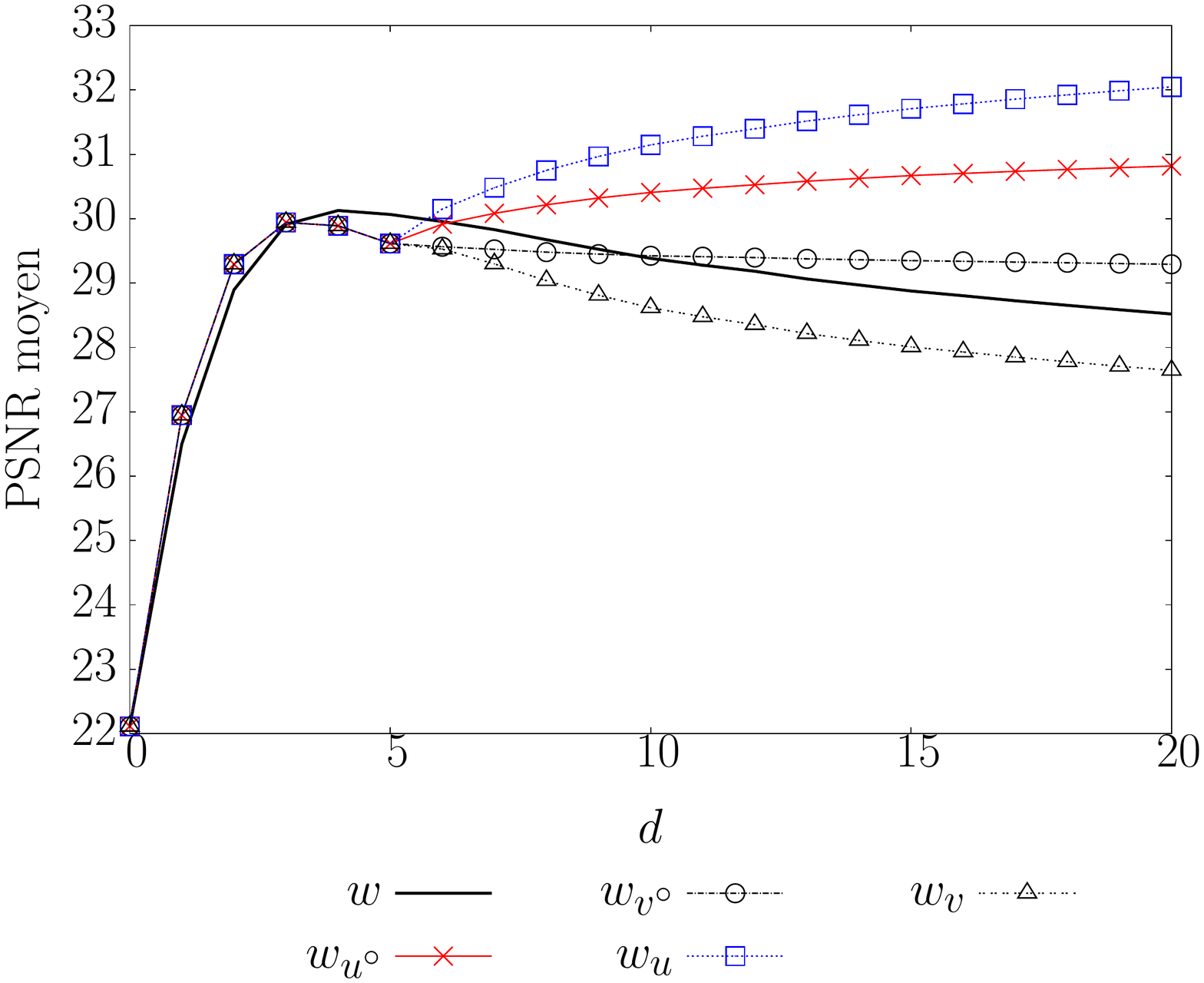}
\end{center}
\legende{PSNR moyen sur 50 images naturelles. Les courbes en rouge (croix) et en bleu (carré) utilisent l'oracle.}
\label{fig:figure1}
\end{figure}

La figure \ref{fig:figure1} donne le PSNR moyen sur 50 images naturelles en fonction du rayon $d$ de la zone de recherche. 
Les images naturelles proviennent du site de traitement d'images en ligne \url{www.ipol.im} (images sous licence {\em Creative Common CC-BY}) 
et des bases de données de Kodak \url{r0k.us/graphics/kodak} et de l'USC-SIPI \url{sipi.usc.edu/database}. 
Elles sont converties en niveaux de gris entre 0 et 255. Nous en affichons quelques-unes à la figure \ref{fig:figure5}.

Les estimateurs $\NLM$ pour les différents poids calculés sans l'oracle, soit $w$, $w_{v}$ et $w_{v^\circ}$, montrent un PSNR moyen en fonction de $d$
qui est une fonction décroissante dès que $d \geq 4$.
Précisons que ce caractère décroissant se retrouve avec d'autres variantes de $\NLM$ non décrites ici, ayant pour objectif d'augmenter avec la distance 
$d$ le nombre de patchs similaires \cite{p2012} (prise en compte du changement de contraste ou de transformations géométriques).
Ces courbes indiquent que l'algorithme de débruitage $\NLM$ est optimal lorsqu'il est restreint à une zone de recherche des similarités très locale. 
Remarquons que le PSNR maximal étant obtenu avec le poids $w$ (pour $d=4$), dans le cadre d'une application effective de débruitage ce poids original 
serait à privilégier sur les variantes calculées par seuillage ou par tri. 

Les courbes correspondantes aux poids obtenus à partir des patchs similaires oracle sont, au contraire, croissantes pour $d \geq 5$, c'est-à-dire dès 
qu'une sélection est effectivement appliquée sur les patchs.
Les poids $\NLM$ calculés sur l'image bruitée sont donc capables d'exploiter les similarités à longue distance, mais à la condition que ces similarités 
soient identifiées sur l'image non bruitée.
Expliquons maintenant cette différence de comportement entre les courbes avec et sans oracle, qui mesure le défaut de non-localité.

\section{Les faux patchs similaires, seuls res\-ponsables du défaut de non-localité}

\begin{figure}[h]
\begin{center}
\includegraphics[width=80mm,bb=4cm 1cm 26cm 19cm,clip=true]{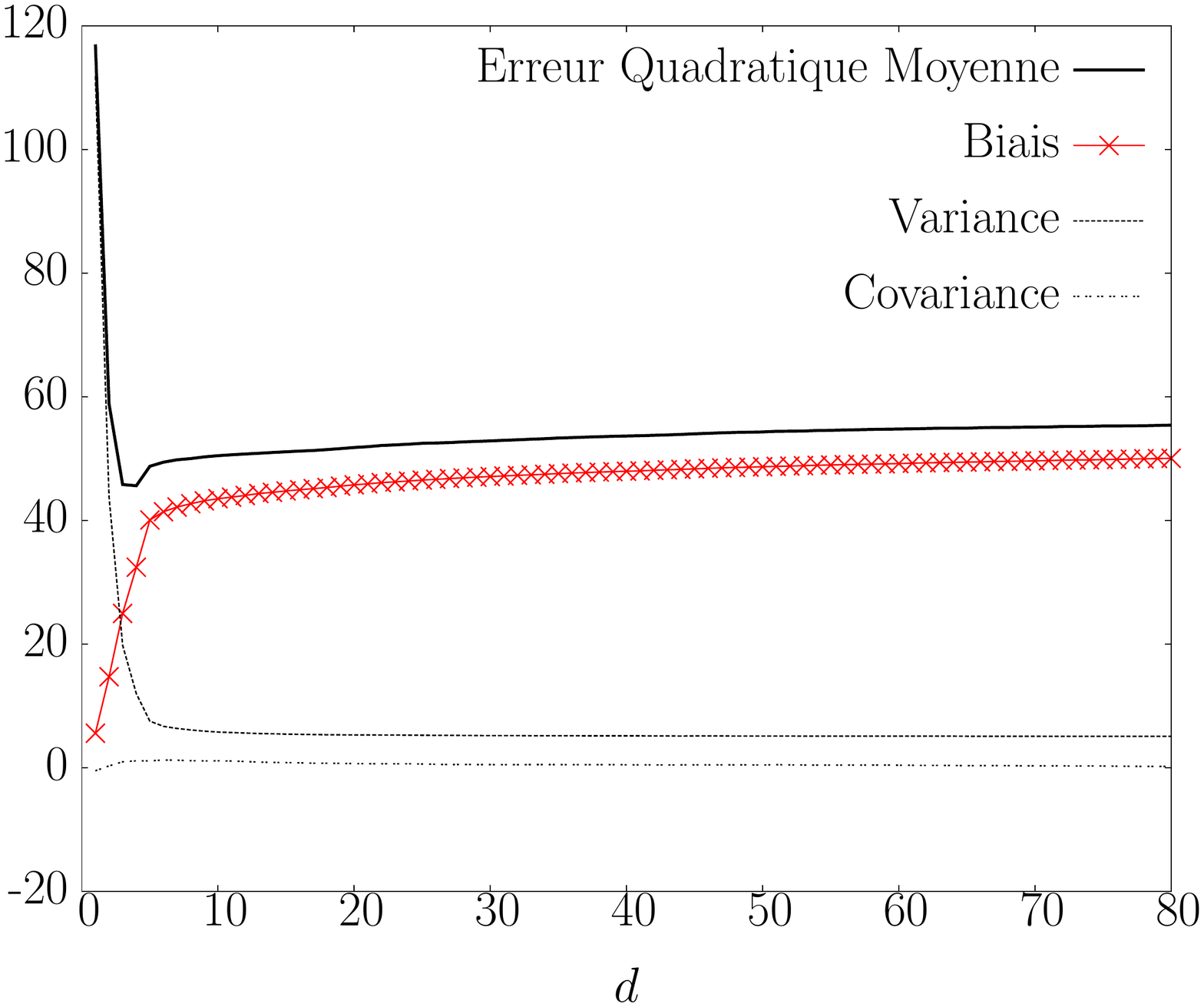}
\end{center}
\legende{Décomposition de l'$\EQM$ sur l'image \textit{Lena}, $f^\circ=v^\circ$: sans oracle le biais est une fonction croissante du paramètre de localité, ce qui entraîne
un minimum local de l'$\EQM$ dans un très petit voisinage (rayon 3 à 4 pixels).}
\label{fig:figure2}
\end{figure}

\begin{figure}[h]
\begin{center}
\includegraphics[width=80mm,bb=4cm 1cm 26cm 19cm,clip=true]{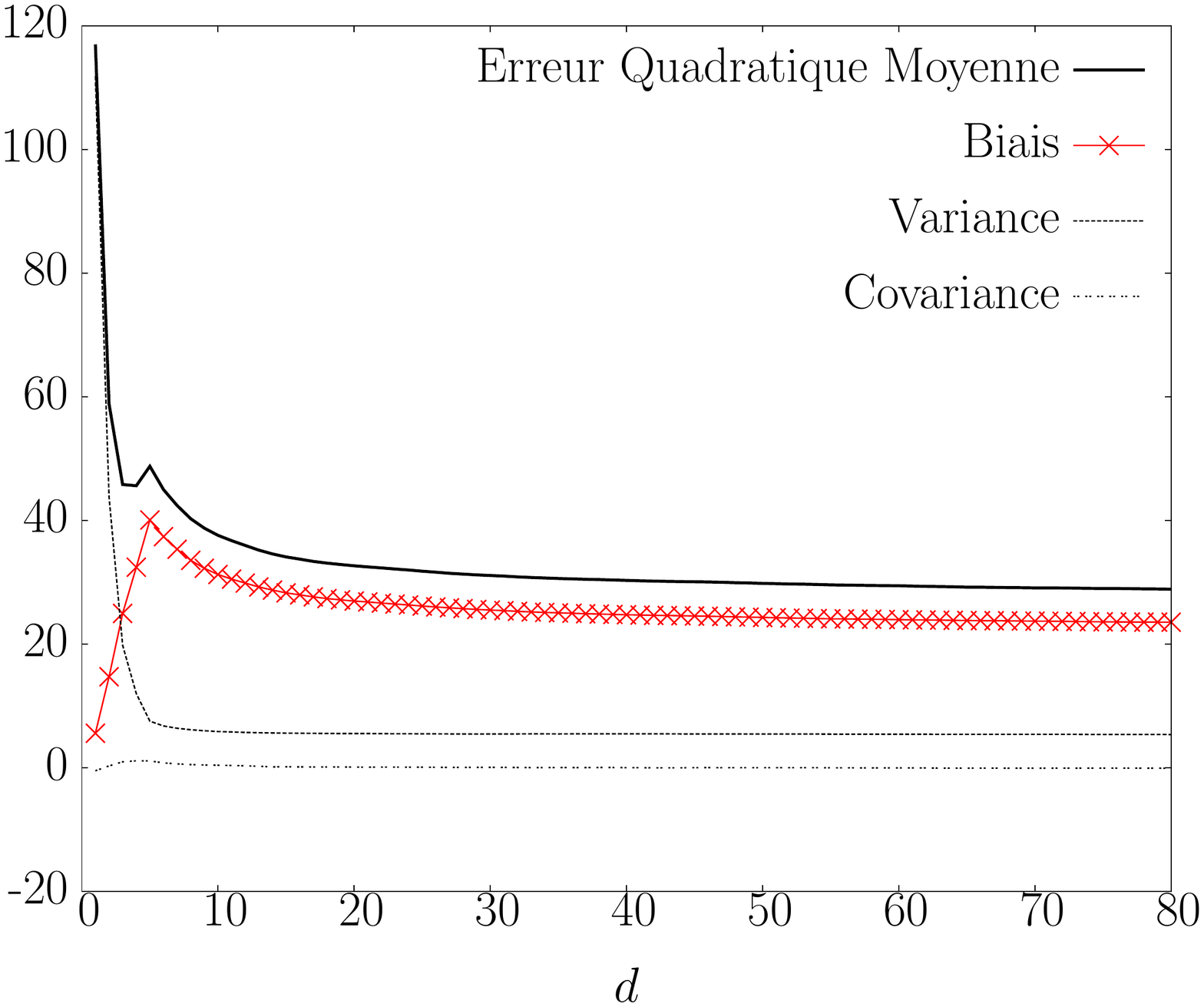}
\end{center}
\legende{Décomposition de l'$\EQM$ sur l'image \textit{Lena}, $f^\circ=u^\circ$: l'oracle donnant les vrais patchs similaires permet au biais de décroître avec la
distance, entraînant la décroissance de l'$\EQM$ (et donc la croissance du PSNR).}
\label{fig:figure3}
\end{figure}

L'Erreur Quadratique Moyenne (EQM) de l'estimateur se décompose suivant les termes $\EQM=\Biais+\Var+\Covar$ où on note $\delta(x,y) = u(y) - u(x)$,
$$\Biais=\frac{1}{|\Omega|}\sum_{x \in \Omega}\left(\sum_{y\in\Oxfronde} w(x,y) \delta(x,y)\right)^2,$$
$$\Var=\frac{1}{|\Omega|}{\sum_{x \in \Omega}\left(\sum_{y\in \Oxfronde} w(x,y) b(y)\right)^2}$$
et le terme $\Covar$ s'exprimant par
$$\frac{2}{|\Omega|}\sum_{x \in \Omega} \left(\sum_{y\in\Oxfronde} w(x,y) b(y)\right) \left(\sum_{y\in \Oxfronde} w(x,y) \delta(x,y)\right).$$
Les figures \ref{fig:figure2} et \ref{fig:figure3} affichent les valeurs de ces trois termes et de leur somme dans le cas $f^\circ=v^\circ$ (b) et $f^\circ=u^\circ$ (c), pour l'image {\em Lena}.
Nous observons que les courbes de variance et de covariance dépendent peu du choix de $f^\circ$; la covariance est quasiment nulle tandis que la variance décroît
rapidement avec la distance $d$ et devient négligeable. 
Seul le biais est responsable du défaut de non-localité, la courbe restant croissante dans le cas $f^\circ=v^\circ$ tandis que dans le cas $f^\circ=u^\circ$ elle décroît
dès que $d \geq 5$.

L'article \cite{bcm2005} précise que le bruit ne perturbe pas en moyenne l'ordre de similarité entre les patchs. En effet,  pour $x \in \Omega$ et $y \in \Omega_x$,
$\mbox{E}\big(\| V(x) - V(y) \|_{2,a}^{2}\big) = \|U(x) - U(y)\|_{2,a}^{2} + 2 \sigma^{2}$.
Mais comme $\|V(x) - V(y)\|_{2,a}^{2} = \|U(x) - U(y)\|_{2,a}^{2} + \|B(x) - B(y)\|_{2,a}^{2} + 2 \big\langle B(x) - B(y), U(x) - U(y) \big\rangle_{2,a}$,
le terme $\|B(x) - B(y)\|_{2,a}^{2} + 2 \big\langle B(x) - B(y), U(x) - U(y) \big\rangle_{2,a}$ perturbe l'ordre de similarité entre les patchs (particulièrement le dernier qui
peut être négatif). Quand le rayon $d$ de la zone de recherche augmente le nombre de patchs augmente, ainsi que le risque de classer comme similaires
des patchs à cause de la seule réalisation du bruit.
L'hypothèse de régularité \og{}$\| V(x) - V(y) \|_{2,a}^{2}$ petit $\Longrightarrow  \|U(x) - U(y)\|_{2,a}^{2}$ petit\fg{} est donc mise en défaut lorsque les deux
patchs sont distants.

\begin{figure}[h]
\begin{center}
\includegraphics[width=80mm,bb=4cm 1cm 26cm 19cm,clip=true]{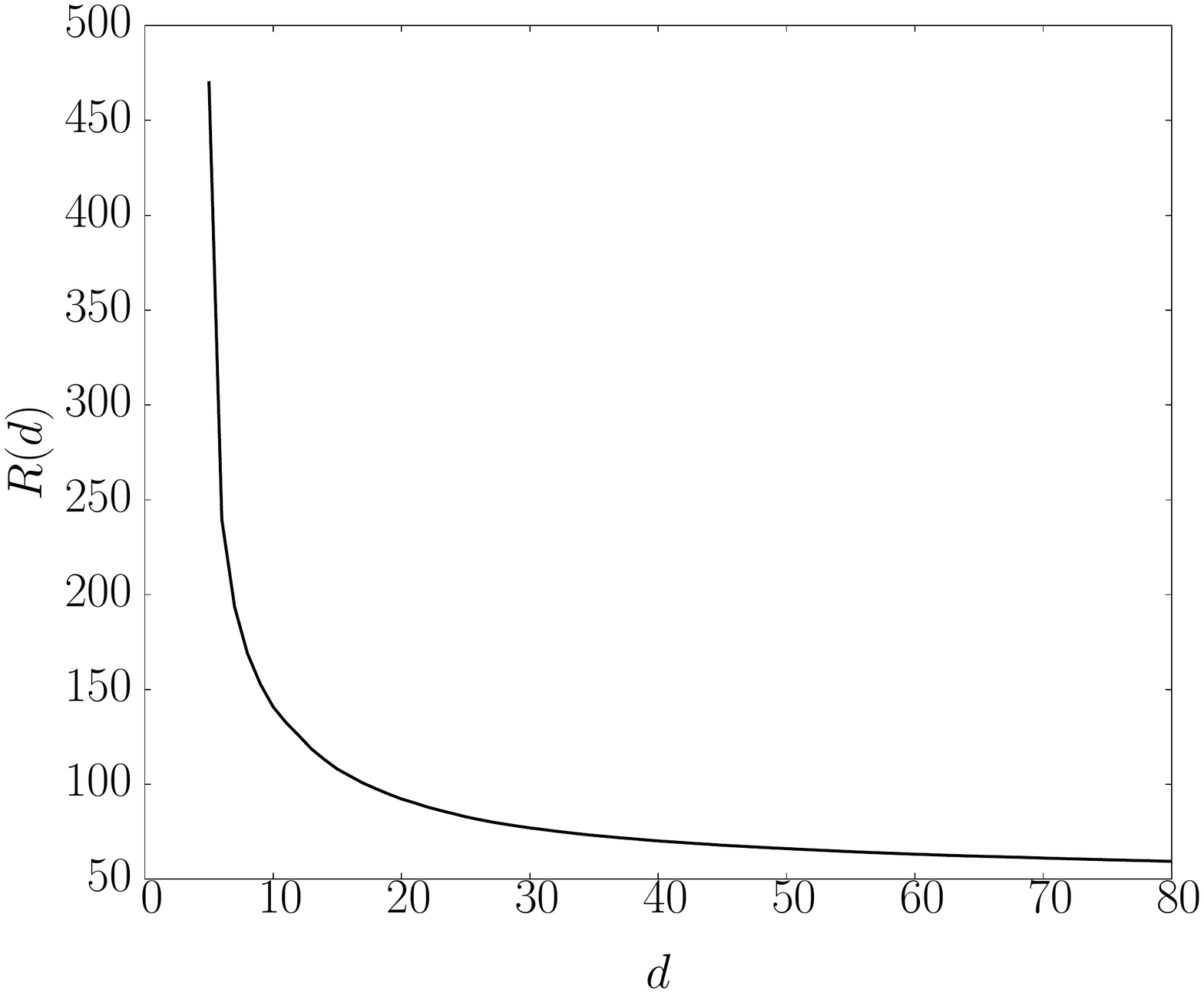}
\end{center}
\legende{Courbe $R(d)$ pour l'image \textit{Lena}, donnant l'erreur quadratique moyenne des pixels considérés comme les plus similaires car au
centre des meilleurs vrais patchs similaires. La décroissance de cette erreur en fonction de la distance établit l'existence d'une régularité
non locale des patchs non bruités aux valeurs centrales.}
\label{fig:figure4}
\end{figure}

\begin{figure}[h]
\begin{center}
\begin{tabular}{cc}
  \includegraphics[width=40mm]{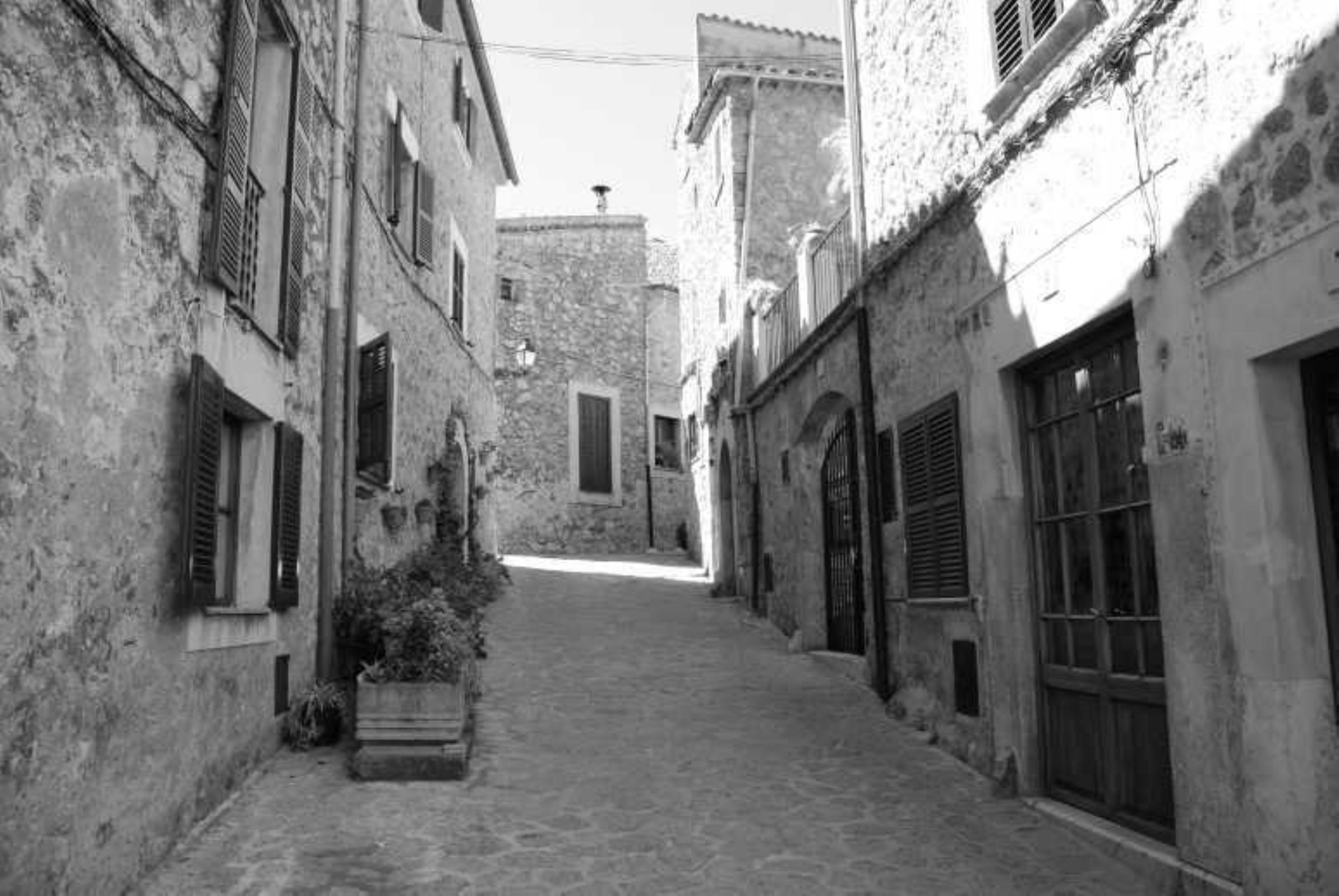}
 & \includegraphics[width=40mm]{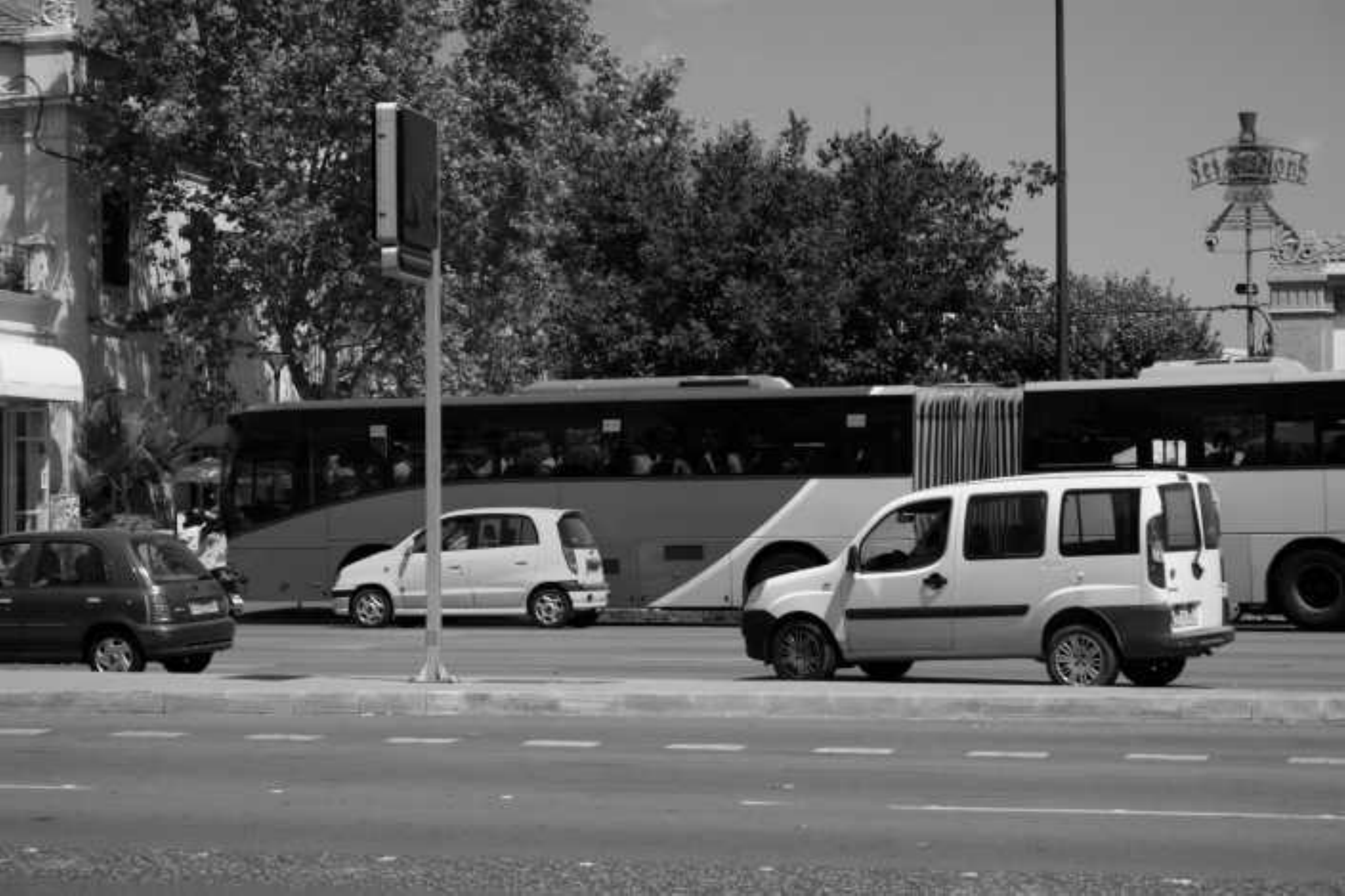}\\
 (a) & (b) \\
\includegraphics[width=40mm]{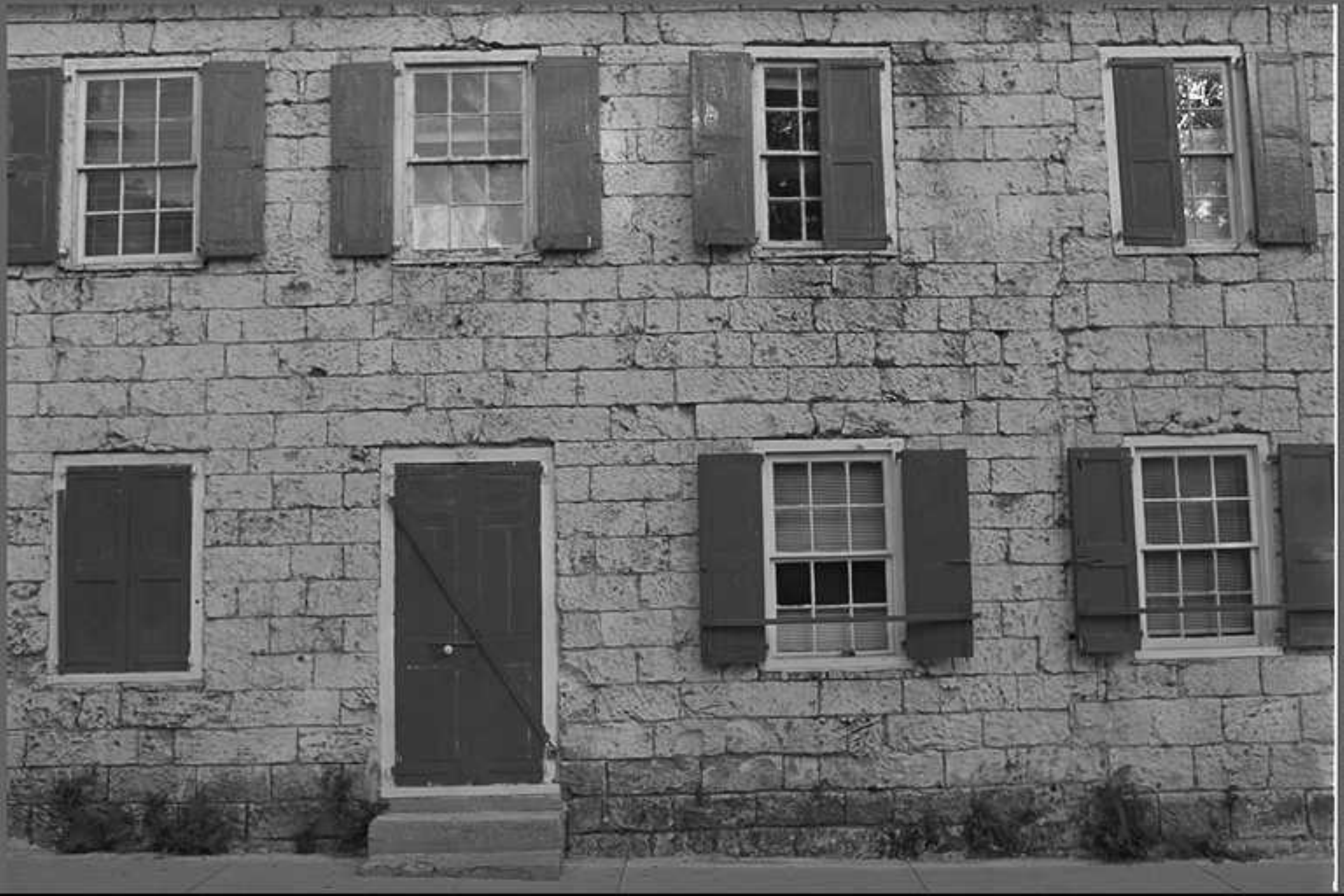}
 & \includegraphics[width=40mm]{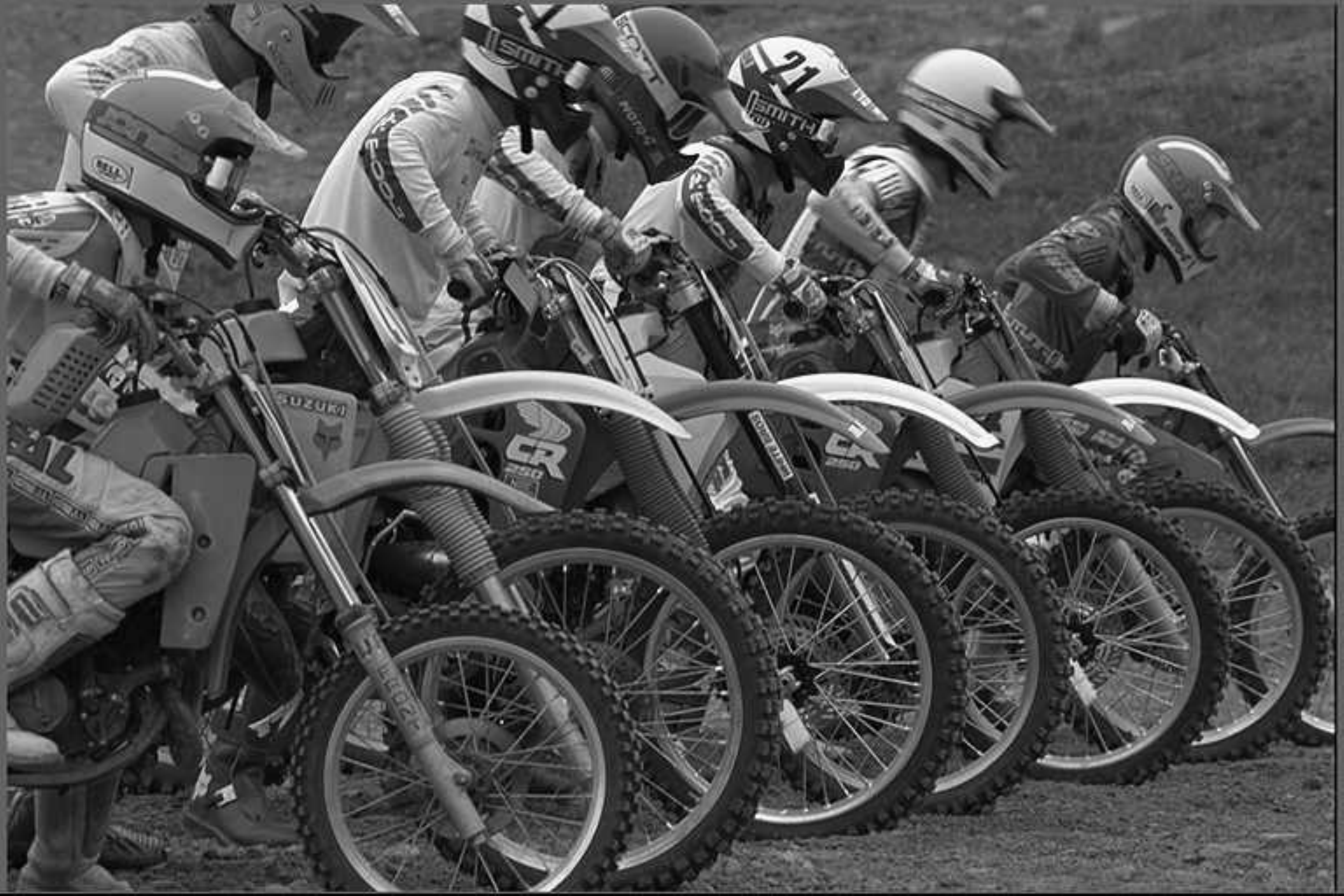}\\
 (c) & (d) \\
  \includegraphics[width=40mm]{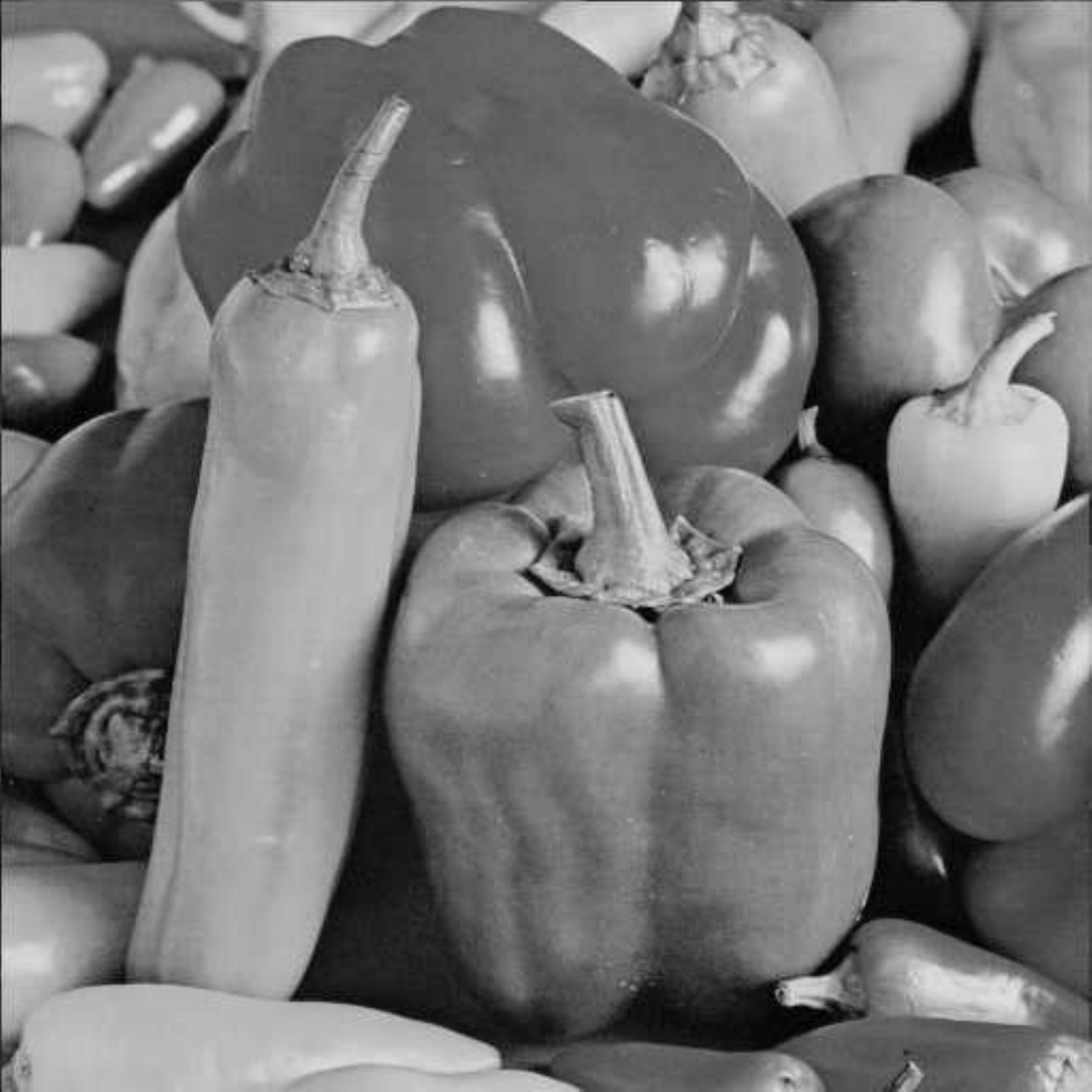}
 & \includegraphics[width=40mm]{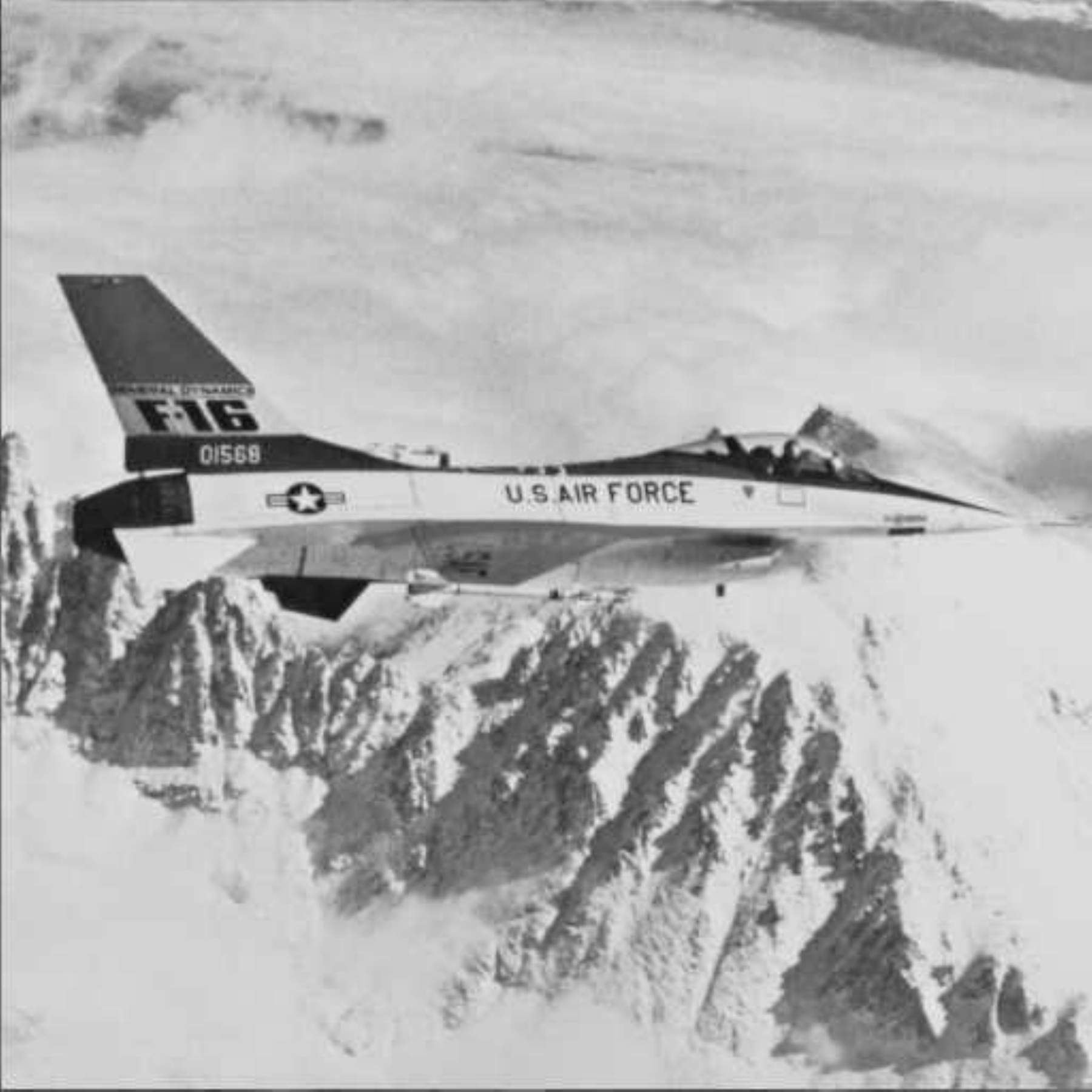}\\
  (e) & (f) \\
   \includegraphics[width=40mm]{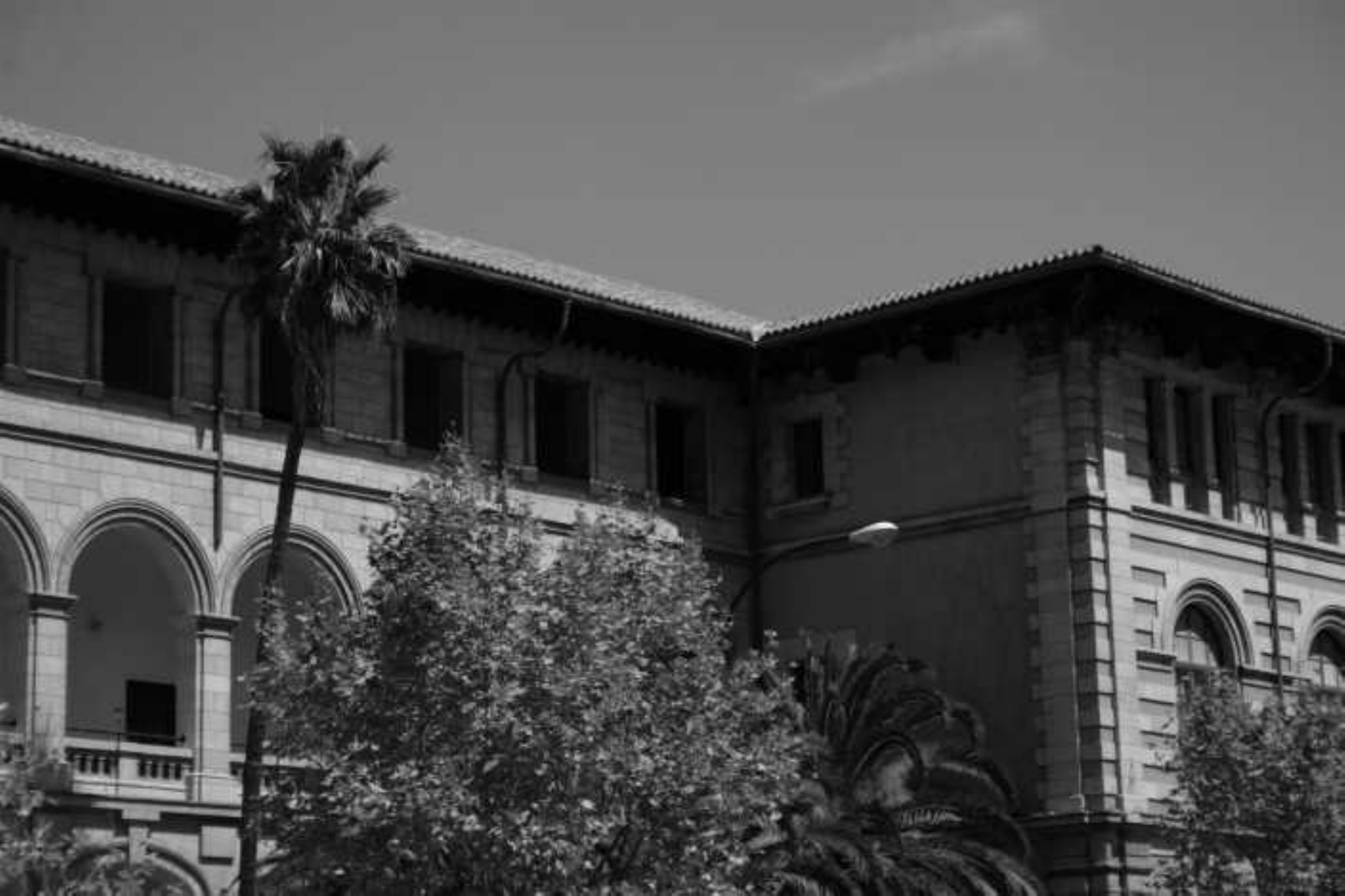}
 & \includegraphics[width=40mm]{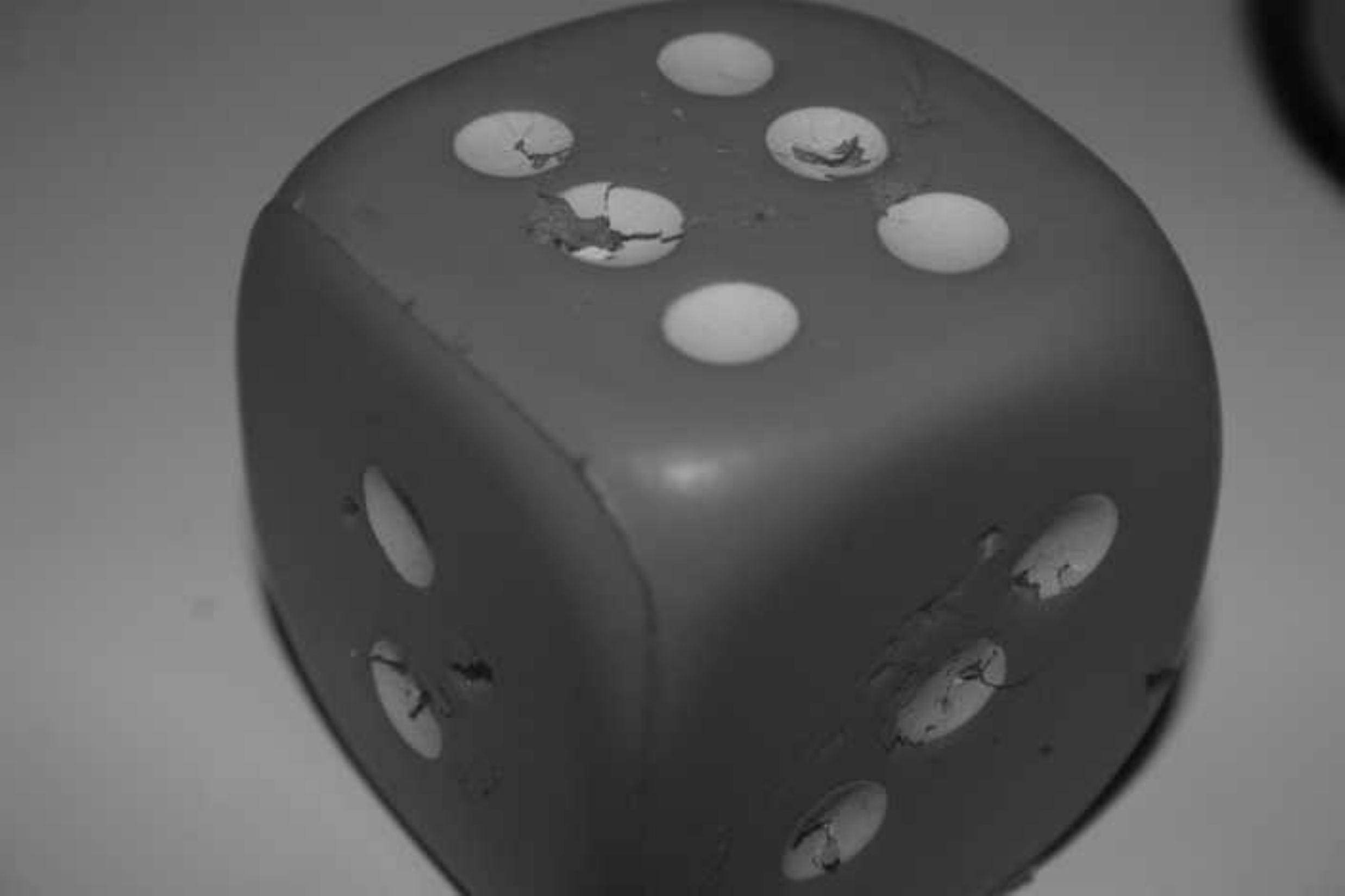}\\
  (g) & (h) \\
    \includegraphics[width=40mm]{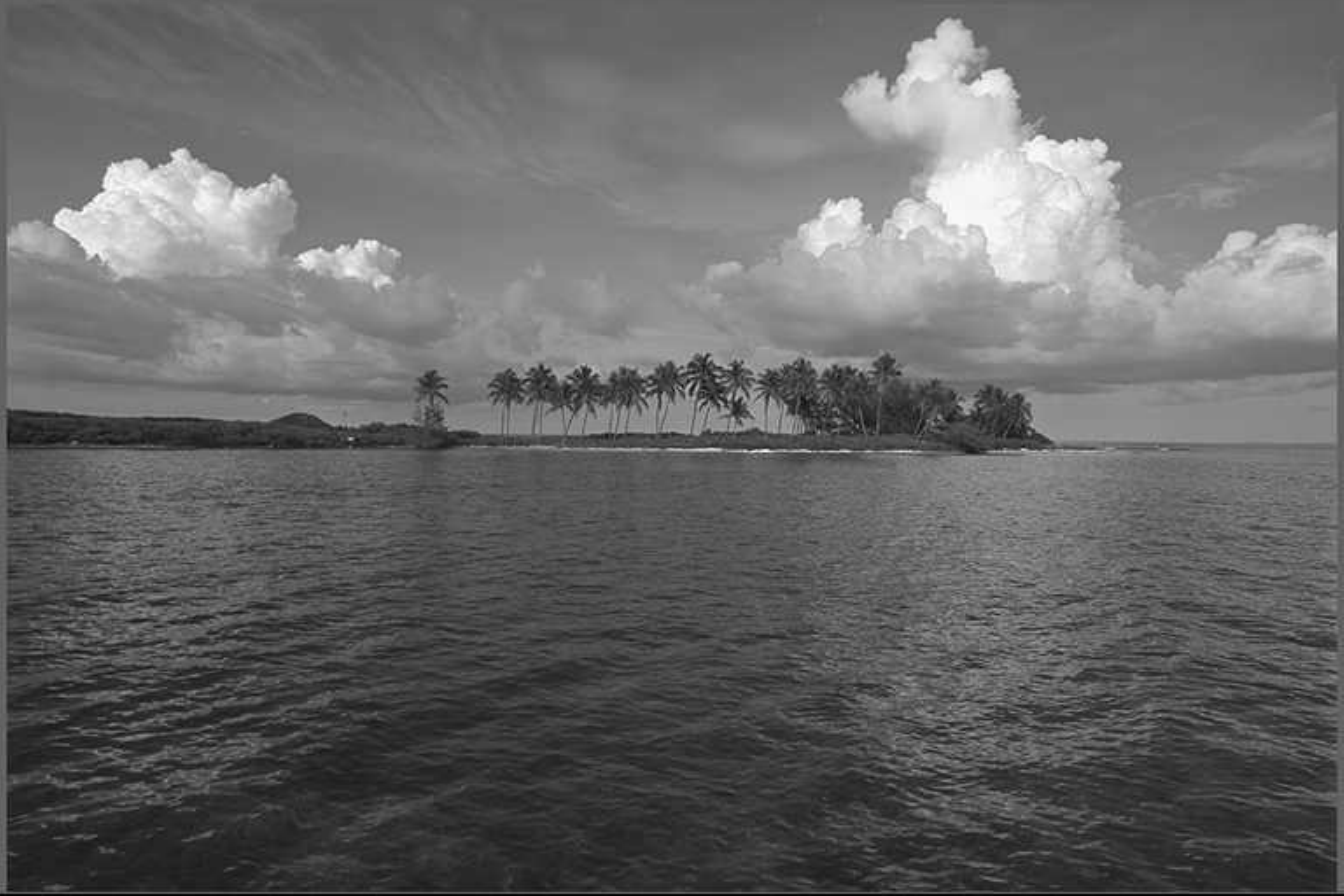}
 & \includegraphics[width=40mm]{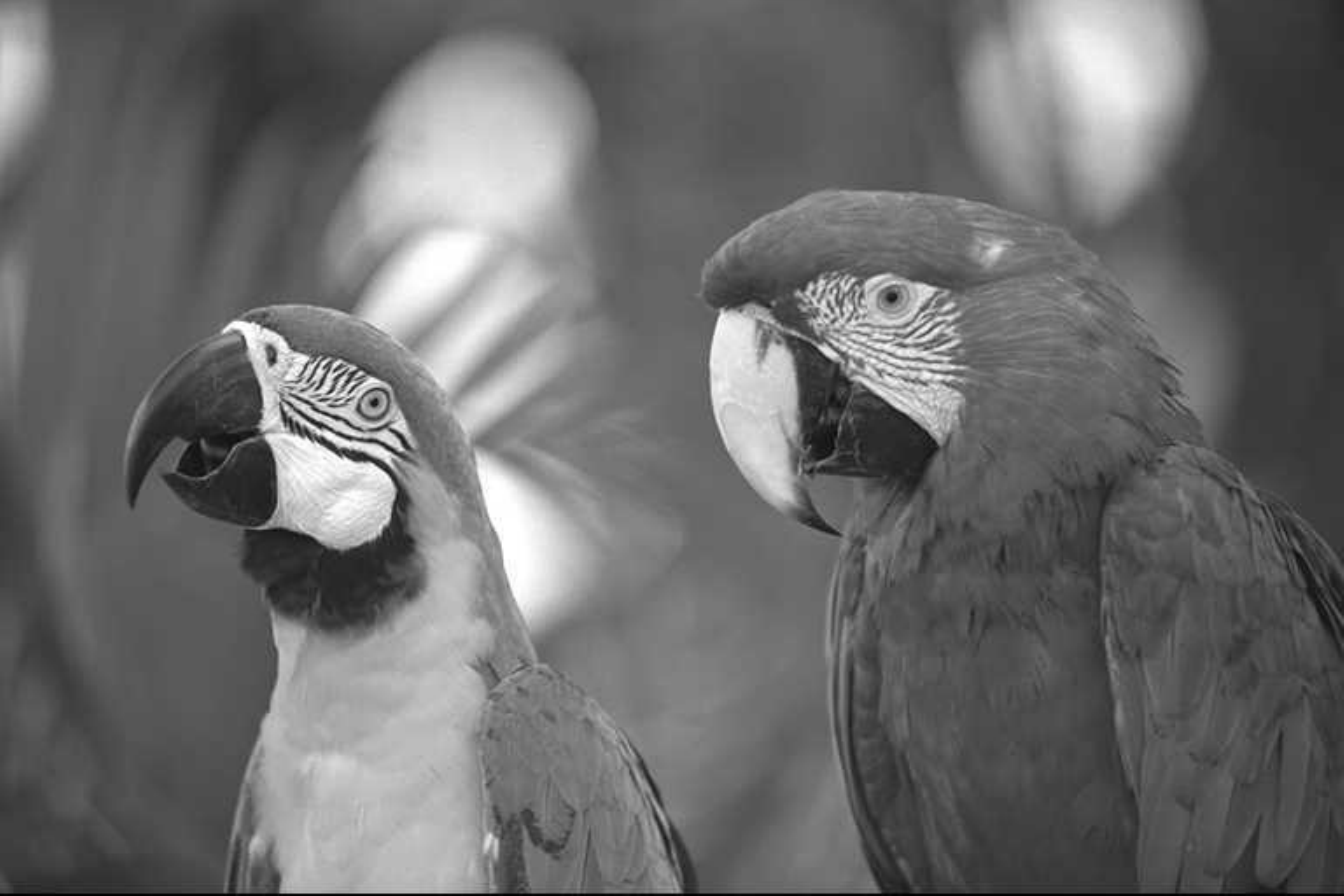}\\
  (i) & (j) \\
\end{tabular}
\legende{Quelques images naturelles considérées comme images originales $u$ pour calculer le PSNR moyen. (a) : CC-BY A. Buades ; (b) : CC-BY M. Colom ;  (c) : Kodak Ref \# JN1206 crédit : Don Cochran ; (d) : Kodak Ref \# R890365 crédit : Steve Kelly ; (e) : Peppers ; (f) : Airplane (F-16) ; (g) : CC-BY M. Colom ; (h) : CC-BY M. Colom ; (i) :  Kodak Ref \# JN0022 crédit : Don Cochran ; (j) : Kodak Ref \# JN1033 crédit : Steve Kelly.}
\label{fig:figure5}
\end{center}
\vspace{-1mm}
\end{figure}

Une autre hypothèse de régularité est implicite dans le modèle $\NLM$ :  \og{}$\|U(x) - U(y)\|_{2,a}^{2}$ petit $\Longrightarrow |u(x)-u(y)|^2$ petit\fg{}. 
La figure \ref{fig:figure4} estime cette régularité en fonction de $d$ à travers le calcul de 
$R(d) = \frac{1}{|\Omega|} \sum_{x \in \Omega}\frac{1}{|\Oxup|} \sum_{y \in \Oxup} |u(x) - u(y)|^2 $.
La fonction étant décroissante, nous concluons que cette propriété est bien vérifiée de manière non locale et ainsi, 
le mauvais choix des patchs similaires est le seul responsable du défaut de non-localité dans la méthode de débruitage $\NLM$.

\section{Conclusion et perspectives}

Le site {\em Google Scholar} référence plus de 3000 textes qui mentionnent l'algorithme $\NLM$; le site {\em IEEE Xplore} dédié aux publications de 
l'institut {\em IEEE} renvoie près de 200 articles comportant cette méthode dans le titre ou dans le résumé. 
La plupart de ces publications considèrent le caractère non local de cet algorithme comme une évidence et cette propriété est parfois évoquée pour 
justifier une nouvelle méthode de traitement des images. Lorsque cette dernière porte sur des données significativement bruitées, notre étude est 
donc susceptible d'invalider le lien de causalité postulé par les auteurs. 
{\em A contrario}, certains articles rejettent l'algorithme de débruitage $\NLM$ au motif de son coût algorithmique considérable lorsqu'il est
effectivement appliqué de manière non locale. 
Cette accusation est infondée: lorsque $\NLM$ est paramétré pour maximiser ses performances, il s'agit d'un algorithme rapide. Aussi, le fait que le paramètre de localité 
$d$ soit généralement choisi plus grand que sa valeur optimale désavantage $\NLM$ dans la comparaison avec d'autres algorithmes de débruitage, y compris avec des variantes de $\NLM$ pour lesquelles d'autres paramètres sont optimisés. 

Au-delà de l'urgence à considérer $\NLM$ pour ce qu'il est réellement (une variante particulièrement performante de filtre à voisinage), se pose la 
question de la capacité d'un algorithme à schéma aussi simple que celui de $\NLM$ à exploiter les similarités non locales des images bruitées. 
Notre étude montre que la distance $\| V(x) - V(y) \|_{2,a}$ n'est pas un critère suffisamment robuste au bruit. Il est tout à fait possible
qu'un critère alternatif permette de reconnaître dans du bruit si deux patchs sont similaires et l'application préalable de celui-ci pour
présélectionner les patchs candidats permettrait alors l'obtention d'un algorithme de débruitage $\NLM$ ``non local''.
Une qualité remarquable de l'algorithme $\NLM$ est de susciter encore, huit ans après son introduction, de telles questions fondamentales.

\end{document}